\documentclass[12pt,leqno]{amsart}
\newtheorem{dummy}{}[section]

\newtheorem{theorem}[dummy]{Theorem}

\newtheorem{lemma}[dummy]{Lemma}

\newtheorem{example}[dummy]{Example}

\begin{document}
\bibliographystyle{plain}
\title{Derived Langlands VIII: Local $GL_{n}$}
\author{Victor P. Snaith}
\date{24 May 2021}
\maketitle
 \tableofcontents 
 
 \section{$GL_{n}$}
 
 Let $K/F$ be a Galois extension of $p$-adic local fields of characteristic zero. Let $\overline{F}$ denote the algebraic closure of $F$. In $GL_{n}K$ let $B_{K}$ denoted the Borel subgroup of upper triangular matrices and let $T_{K}$
 be the torus of diagonal matrices. Let $L$ be the rational maps on $T_{K}$ and if $\lambda \in L$ we have
 \[ \lambda : \left(  \begin{array}{cccccc}
 t_{1} & 0 & 0 & \ldots & \ldots & 0 \\
 0 &  t_{2} & 0 & \ldots & \ldots & 0 \\
 0&  0 &  t_{3} &  \ldots & \ldots & 0 \\
 0 & \ldots & \ldots & \ldots & \ldots &0 \\
  0 & \ldots & \ldots & \ldots & \ldots &0 \\
   0 & 0 &0& \ldots & \ldots & t_{n}  \\
 \end{array} \right)  \longrightarrow (t_{1}^{m} \ldots t_{n}^{m} ) \]
 with $m_{1}, \ldots , m_{n} \in {\mathbb Z}$ we write $\lambda = (m_{1}, \ldots , m_{n})$ so $L \cong  {\mathbb Z}^{n}$.
 
 The Lie algebra of $GL_{n}$ consists of the $n \times n$ matrices and the roots are
 \[ \begin{array}{ll}
 a_{1} &= (1, -1, 0, \ldots, \ldots, 0) \\
  a_{2} &= (0, 1, -1, 0, \ldots,  0) \\
  & \vdots \hspace{15pt} \hspace{15pt}   \vdots  \hspace{15pt}  \vdots \\
  a_{n} & = (0,0, \ldots , 0 , 1 , -1). \\
 \end{array} \]
 
 Let $(-,-)$ be the familiar bilinear form on the $n \times n$ matrices with entries in $K$.
 
 \section{Preamble}
 
 I shall describe how the material of this note constructs a 2-variable L-function for two representations of general linear groups of a $p$-adic local field $F$ defined over $\overline{F}$, the algebraic closure of $F$ and admissible in the sense of \cite{Sn20}. More precisely, admissibility in the above sense has the following meaning. 
 We shall say that $V$ is ${\mathcal M}_{cmc, \underline{\phi}}(G)$-smooth\footnote{In previous versions I made a goofy typographical blunder and wrote $\bigcup$ instead of ${\rm Span}$.} if 
 \[ V = {\rm Span}_{ (H, \phi) \in {\mathcal M}_{cmc, \underline{\phi}}(G)} \ V^{(H, \phi)}.\]
 In addition we shall say that $V$ is ${\mathcal M}_{cmc, \underline{\phi}}(G)$-admissible if 
 ${\rm dim}_{k}V^{(H, \phi)} < \infty$ for all  $(H, \phi) \in {\mathcal M}_{cmc, \underline{\phi}}(G)$.
 
 I apologise for the state of these notes but I have health problems \cite{JMOS18}  and impaired vision so this is the best I can do.
 
 We begin with representations $\rho$ and $\pi$ of $GL_{n}F$ and $GL_{m}F$ respectively over $\overline{F}$ which are admissible in the above sense.
 
 We choose an unramified Galois extension of local fields $K/F$ and let $\tilde{\pi}$ be the admissible representation of the semi-direct product ${\rm Gal}(K/F) \propto GL_{m}K$ given by base-change \cite{AC89}.
 
 The Galois semi-direct product acts simplicially on the Bruhat-Tits building, for example with the simplicial structure described in \cite{PG97}. As illustrated in \S 4, \S 5, \S 6 the fundamental subcomplex is a finite number of simplices whose stabiliser is $K^{*} U_{K} B_{K}$ where $B_{K}$ is a Borel subgroup, $U_{K}$ is compact open and $U_{K} B_{K}$ satisfies the cohomology of \S 7 and Lemma \ref{5.1a}.
 
 The observation of Digne-Michel in \S 4 matches up  conjugacy classes  in  $GL_{m}F$ with conjugacy classes  of elements $( \rm{  Frob}_{F}, A)$ in the Galois semi-direct product. Therefore the object of the exercise is to construct an L-function  $ L(s, \rho, \tilde{\pi})$ depending on these conjugacy classes only. A multiplicative Euler characteristic of the form  $ L(s, \rho, \tilde{\pi})$ may then give a well-defined  $ L(s, \rho, \pi)$.

This paper explains why $L(s,\rho, \pi)$ in my di-$p$-adic set up should be given by Artin's usual determinantal formula (see \S 9 {\bf Part B}).

 \section{ Lie algebras and roots} 
  \vspace{10pt} 
 
 A generalised Cartan matrix is a square matrix $A = a_{i,j}$ with integer entries in $\Sigma_{n} \int D$ where $D$ is diagonal matrices such that (i) each $a_{i,i} = 2$ , \ (ii) if $i \not=j$ then $a_{i,j} \leq  0$ \ (iii) $a_{i,j} = 0$ if and only if $a_{j,i}=0$. 
 
 For example, for the group $G_{2}$ the Cartan matrix is 
 \[  \left(  \begin{array}{cc}
 2 & -3 \\
 \\
 -1 & 2 
 \end{array}   \right) = \left(  \begin{array}{cc}
 3 & 0 \\
 \\
 0 & 1 
 \end{array}   \right)  \left(  \begin{array}{cc}
 \frac{2}{3} & -1 \\
 \\
 -1 & 2 
 \end{array}   \right)  = DS \]
 We can always choose $D$ to have positive diagonal entries. In that case, if $S$ is positive definite then $A$ is said to be a Cartan matrix.
 
 The Cartan matrix of a simple Lie algebra is the matrix whose entries are the scalar products
 \[   a_{j,i} = 2 \frac{(r_{i}, r_{j}) }{(r_{j}, r_{j}) }  ,  \]
 sometimes called the Cartan integers, where the $r_{i}$'s are the {\bf simple roots} of the algebra. Condition (i) is obvious. Condition (ii) follows from the fact that  if $i \not= j$ then $r_{j} - 2 \frac{(r_{i}, r_{j}) }{(r_{j}, r_{j}) } r_{i} $ is a linear combination of the simple roots $r_{i}$ and $r_{j}$ with a positive coefficient for $r_{j}$ and so the coefficient of $r_{i}$ has to be non-negative. Condition (iii) follows because orthogonality is a symmetric relation. Finally
 \[ D_{i,j} = \frac{ \delta_{i,j}}{ (r_{i}, r_{i})  } , \ S_{i,j} = 2(r_{i}, r_{j}). \]
 $S$ is positive definite because the simple roots span a Euclidean space.
 
 Conversely, from the Cartan matrix, one can reconstruct the Lie algebra.
 
 Recall that a Lie aqlgebra ${\mathcal G}$ is a vector space with a Lie product
 $(x,y) \mapsto [x,y]$ which satisfies $[x,y] = - [y,x]$ and the Jacobi identity 
 \[   [x, [y, z]] +  [z, [x,y]] + [y , [z,x]] = 0.\]
 
 Let $E$ be a finite-dimensional Euclidean vector space, with standard inner product $(-,-)$. A root system $\Phi$ in $E$ is a set of non-zero vectors (i) spanning $E$, \ (ii) The only scalar multiples of $\alpha \in \Phi$ are $\pm \alpha$, \ (iii) for every root $\alpha \in \Phi$ the set $\Phi$ is closed under reflection in the hyperplane perpendicular to $\alpha$, \ (iv) if $\alpha, \beta \in \Phi$ then the projection of $\beta$ onto the line through $\alpha$ is an integer multiple or a half-integer multiple of $\alpha$.
 
 Equivalently to (iii) and (iv) are the conditions
 \[ (iii)' \  {\rm if} \   \alpha, \beta \in \Phi  \ {\rm then}  \ \sigma_{\alpha}(\beta) =
  \beta - 2 \frac{( \alpha , \beta )}{( \alpha ,  \alpha )} \alpha  \in \Phi \]
 and 
 \[ (iv)' \  {\rm if} \   \alpha, \beta \in \Phi  \ {\rm then} \   \langle  \beta , \alpha \rangle = 2  \frac{( \alpha , \beta )}{( \alpha ,  \alpha )} \in {\mathbb Z}. \]
 
 Some authors only include conditions (i)-(iii) in the definition of a root system. In this context a root system which also satisfies (iv) is known as crystallographic and the systems satisfying (ii) are called reduced. Here all root systems will be both crystallographic and and reduced.
 
 Note that the pairing $ (\alpha, \beta) \mapsto  \langle  \beta , \alpha \rangle$ - $ \Phi \times \Phi \longrightarrow {\mathbb Z}$ - is {\em not} an inner product because it is not symmetric in the variables and is linear only in $\beta$.

\section{The observation of Digne-Michel \cite{DM85}, \cite{DiMi} and (\cite{Sn18} pp.267-269)}

\begin{dummy}
\label{I10.1}
\begin{em}

Let $F$ be a $p$-adic local field of characteristic zero with algebraic closure $\overline{F}$ and within that a maximal unramified extension inside which I shall assume all our important field elements lie.
Let ${\rm Frob}_{F}$ be the Frobenius automorphism topologically generating the absolute 
unramified Galois group of $F$,
 ${\rm Gal}( F_{{\rm ur}}/ F)$, so that
${\rm Frob}_{F}^{n}$ generates ${\rm Gal}( F_{{\rm ur}}/ K)$ with $[K:F]=n$. 

Assuming that the Lang-Steinberg Theorem holds for $U_{F}$ which is a subspace of $GL_{s}K$ acted upon by the unramified Galois group, for any $Y \in GL_{s} {F}$ there exists an unramified extension $K$ and a matrix $X \in GL_{s}K$ such that $Y = X^{-1} {\rm Frob}_{F}(X)$. 

Observe that, if $[K:F]=n$,

\[ \begin{array}{l}
  X^{-1} {\rm Frob}_{F}(X)  \in  GL_{s}K  \\
\\
 \Longleftrightarrow   \   {\rm Frob}_{F}^{n}(  X^{-1}) {\rm Frob}_{F}^{n+1}(X) =  X^{-1} {\rm Frob}_{F}(X)  \\
\\
  \Longleftrightarrow  X {\rm Frob}_{F}^{n}(  X^{-1})  =   {\rm Frob}_{F}(X  {\rm Frob}_{F}^{n}(  X^{-1}) )  \\
\\
 \Longleftrightarrow     X {\rm Frob}_{F}^{n}(  X^{-1})  \in GL_{s}F.
\end{array} \]

Let $W \in GL_{s}$ be another matrix and suppose that 
\[  (  {\rm Frob}_{F} , X^{-1}  {\rm Frob}_{F}(X)) ,  \   (  {\rm Frob}_{F} ,  W^{-1}  {\rm Frob}_{F}(W))   
\in {\rm Gal}( K /F )  \propto GL_{s}K \]
are conjugate by $(1 , V) \in GL_{s}K$. Hence, say,
\[ \begin{array}{l}
 (1 , V)  ( {\rm Frob}_{F} , X^{-1}  {\rm Frob}_{F}(X)) (1, V^{-1})  \\
\\
=  (  {\rm Frob}_{F} , V X^{-1}  {\rm Frob}_{F}(X)) (1, V^{-1})  \\
\\
 =  ( {\rm Frob}_{F}, V X^{-1}  {\rm Frob}_{F}(X V^{-1}) )  \\
\\
 =  (  {\rm Frob}_{F} , W^{-1}  {\rm Frob}_{F}(W)).
\end{array} \]

Therefore $WVX^{-1} =  {\rm Frob}_{F}( WVX^{-1} ) \in GL_{s}F$ and
\[  \begin{array}{l}
 (1 , WVX^{-1}) (  {\rm Frob}_{F}^{n} ,   X {\rm Frob}_{F}^{n}(  X^{-1})) ( 1 , XV^{-1}W^{-1})  \\
\\
= (  {\rm Frob}_{F}^{n} ,   WVX^{-1} X  {\rm Frob}_{F}^{n}(  X^{-1})) ( 1 , XV^{-1}W^{-1}) \\
\\
=  (  {\rm Frob}_{F}^{n} ,   WVX^{-1} X  {\rm Frob}_{F}^{n}(  X^{-1} XV^{-1}W^{-1})) \\
\\
= ( {\rm Frob}_{F}^{n} ,   WV  {\rm Frob}_{F}^{n}( V^{-1}W^{-1})) \\
\\
= (  {\rm Frob}_{F}^{n} ,   WV V^{-1} {\rm Frob}_{F}^{n}( W^{-1})) \\
\\
=  ( {\rm Frob}_{F}^{n} ,   W  {\rm Frob}_{F}^{n}( W^{-1})) .
\end{array} \]
Therefore 
\[ (  {\rm Frob}_{F}^{n} ,   X  {\rm Frob}_{F}^{n}(  X^{-1})),   \   (  {\rm Frob}_{F}^{n} ,   W  {\rm Frob}_{F}^{n}( W^{-1})) \in
 {\rm Gal}( K / F )  \times   GL_{s}K  \]
are conjugate by an element of $GL_{s} F$ or, equivalently, 
\[   X  {\rm Frob}_{F}^{n}(  X^{-1}),   \     W  {\rm Frob}_{F}^{n}( W^{-1}) \in
   GL_{s}F \]
are conjugate in $ GL_{s}F$.

Conversely, if there exists $A \in GL_{s}F$ such that
\[  X \Sigma^{n}(  X^{-1}) = A  W  {\rm Frob}_{F}^{n}( W^{-1}) A^{-1} \in GL_{s}F \]
then
\[ W^{-1} A^{-1} X =  {\rm Frob}_{F}^{n}( W^{-1} A^{-1} X ) \in GL_{s}K.  \]
Also
\[ \begin{array}{l}
  (1 ,  W^{-1} A^{-1} X ) (  {\rm Frob}_{F} , X^{-1}  {\rm Frob}_{F}(X)) (1 , X^{-1}AW ) \\
\\
 =   (  {\rm Frob}_{F} ,  W^{-1} A^{-1} X X^{-1}  {\rm Frob}_{F}(X)) (1 , X^{-1}AW ) \\
\\
=  (  {\rm Frob}_{F} ,  W^{-1} A^{-1} X X^{-1}  {\rm Frob}_{F}(X X^{-1}AW) ) \\
\\
= (  {\rm Frob}_{F},  W^{-1} A^{-1}  {\rm Frob}_{F}(AW) ) \\
\\
= (  {\rm Frob}_{F} ,  W^{-1} A^{-1}A  {\rm Frob}_{F}(W) )  \\
\\
= (  {\rm Frob}_{F} ,  W^{-1}   {\rm Frob}_{F}(W) )
\end{array} \]
so that 
\[  (  {\rm Frob}_{F} , X^{-1}  {\rm Frob}_{F}(X))  ,  \   ( {\rm Frob}_{F} ,  W^{-1}   {\rm Frob}_{F}(W) ) 
\in {\rm Gal}( K / F)  \propto GL_{s} K \]
are conjugate by an element of $GL_{s} K $.

Therefore we have proved the following result.
\end{em}
\end{dummy}
\begin{theorem}{$_{}$}
\label{I10.2}
\begin{em}

There is a one-one correspondence of the form
\[  \left\{  \begin{array}{c}
  {\rm conjugacy} \\
\\
 {\rm classes \ in}  \\
\\
 GL_{s} F
\end{array} \right\}   \leftrightarrow  
\left\{  \begin{array}{c}
 GL_{s}K - {\rm conjugacy} \\
 \\
 {\rm  classes \ of}  \\
\\
{\rm elements}  \  (  {\rm Frob}_{F} , A) \ {\rm in} \\
  \\
{\rm Gal}( K / F )  \propto  GL_{s} K 
\end{array}  \right\} \]
given by
\[  X  {\rm Frob}_{F}^{n}(X^{-1})  \leftrightarrow  (  {\rm Frob}_{F} , X^{-1}  {\rm Frob}_{F}(X))  \]
for $X \in GL_{s}M$ for some $M/F$ with $[M:F]=t$.
\end{em}
\end{theorem}

\section{$\pi$ on $GL_{1}$}
 
 In (\cite{LLProb1970} Lemma 1, p.30), $G = GL_{1}$, $G_{F} = F^{*}$, $G_{K} = K^{*}$ we have $K^{*} = G_{K} = B_{K}U_{K}$, $U_{K} = {\mathcal O}_{K}^{*}$\footnote{$U_{F}$ must be compact if $f \not= 0$ and $f$ is compactly supported.}.
 and the vanishing of $H^{1}({\rm Gal}(K/F); U_{K})$ and $H^{1}({\rm Gal}(K/F); B_{K} \bigcap U_{K})$
 implies that $B_{K} = K^{*} = U_{K}$ and $B_{F} = F^{*},  U_{F} = {\mathcal O}_{F}^{*}$. The set of roots is empty for $GL_{1}$. $T^{K} = K^{*}, T_{F}=F^{*}$.
 
 $L = \{ x \in F^{*} \ | \ x \mapsto x^{n} \} \cong {\mathbb Z}$, rational functions (see \cite{LLProb1970}, p.33).
 $\hat{L} = {\rm Hom}(L , {\mathbb Z}) \cong {\mathbb Z}$,  ((see \cite{LLProb1970}, pp.34/35).
 Also $\hat{M} = \hat{L}^{{\rm Gal}(K/F) }$ and $\Lambda(\hat{M})$ is the $k$-group ring of $\hat{M}$. 
  
  $G_{K} = B_{K} = K^{*}$ and $v = v_{F} : F^{*} \longrightarrow {\mathbb Z}$ is the normalised valuation on $F$.
  $U_{K} = {\mathcal O}_{K}^{*}$ which is compact, as required. ${\mathcal O}_{K}/{\mathcal P}_{K} = {\mathbb F}_{q^{d}}$ with residue degree $[  {\mathbb F}_{q^{d}} :  {\mathbb F}_{q} ] = d$. $G_{F} = B_{F} = F^{*}$ and
  ${\mathcal O}_{F}/{\mathcal P}_{F} = {\mathbb F}_{q}$. We have $\rho : F^{*} \longrightarrow  \overline{F}^{*}$ and $ \pi : K^{*}  \longrightarrow  \overline{F}^{*}$.

  We have $({\mathcal O}_{K}^{*}, \phi)^{K^{*}}$ so by analogy the the case $\phi =1$ in \cite{LLProb1970} and similarly with $K$ replaced by $F$. We define, for $K/F$ large enough
 \[   \begin{array}{l}
 C_{c}(K^{*},  ({\mathcal O}_{K}^{*}, \phi)^{K^{*}}) = \{ f : K^{*} \longrightarrow \overline{F}  \  |  \ f \ {\rm compactly \ supported} \\
 \hspace{60pt} f(ug) = f(gu) = \phi(u)^{-1} f(g) , \ u \in {\mathcal O}_{K}^{*}, g \in K^{*} \}  .
 \end{array} \]
 and similarly with $K$ replaced by $F$. 
 
 Let $v$ be the normalised $F$-adic valuation on any extension $K/F$ so if we have $\lambda = (m)$ in $\hat{L} = {\rm Hom}(L, {\mathbb Z})$ with $t \in F$ we have $\lambda(t) = t^{m}$ and if $v(t) = s$ then $t = \pi_{F}^{s} u$ with $u \in U_{F} = {\mathcal O}_{F}^{*}$. 
 
 We have 
 \[  C_{c}(K^{*},  ({\mathcal O}_{K}^{*}, \phi)^{K^{*}})  \longrightarrow  \overline{F}  \]
 the details of which will appear later.
 \newpage

\section{$\pi$ on $GL_{2}$}

The Galois action on $K \oplus K$ sends stabilisers of homothety types to themselves. Hence ${\rm Gal}(K/F)$ acts on stabilisers of simplices in the building. I believe this is true in general for the action on the Tammo Tom Dieck space of $GL_{n}K$. The monomial resolution of an admissible representation of $GL_{2}$ is given by the double complex of putting the bar resolution of the hyperHecke algebra of the restriction on $\pi$ to a compact open mod the centre stabiliser of a simplex. The fundamental domain of the $GL_{2}$-action on the building is a $1$-simplex which has $2$ vertices with compact open modulo the centre stabilisers. The first vertex has stabiliser $G_{K} = K^{*} \cdot GL_{2}{\mathcal O}_{K}$ whose elements look like
\[ G_{K} = \{ \pi_{K}^{i}   \left(  \begin{array}{cc}
a & b \\
\\
c & d \\
\end{array} \right)  \ | \ ad - bc \in {\mathcal O}_{K}^{*}, a,b,c,d \in  {\mathcal O}_{K}, \ i \in {\mathbb Z}  \}  .  \]

We have a Chevalley subgroup $U_{K}$ as in (\cite{LLProb1970} Lemma 1, p.30). Inside $G_{K}$ take the abelian group
\[    U_{K} = \{ \left(  \begin{array}{cc}
\pi_{K}^{i} a & 0 \\
\\
\pi_{K}^{i} c & \pi_{K}^{i} a
\end{array} \right) \ | \ a \in  {\mathcal O}_{K}^{*}, c \in {\mathcal O}_{K}  \}      \]
and the Borel subgroup
\[    B_{K} = \{ \left(  \begin{array}{cc}
\pi_{K}^{j} \alpha   & \pi_{K}^{j} \beta   \\
\\
0 & \pi_{K}^{j} \delta  \end{array} \right)  \ | \   \alpha, \delta \in  {\mathcal O}_{K}^{*}, \beta \in {\mathcal O}_{K}   \}      \]
so that $G_{K} = U_{K} B_{K}$ and $U_{K}$ is self-normalising in $G_{K}$ as the following calculation shows.
\[  \begin{array}{l}
\left(  \begin{array}{cc}
\pi_{K}^{i} a & 0 \\
\\
\pi_{K}^{i} c & \pi_{K}^{i} a
\end{array} \right)  \left(  \begin{array}{cc}
\pi_{K}^{j} \alpha   & \pi_{K}^{j} \beta   \\
\\
0 & \pi_{K}^{j} \delta  \end{array} \right)  = \left(  \begin{array}{cc}
a  \alpha   & a \beta   \\
\\
c \alpha  & c \beta + a  \delta  \end{array} \right) \pi_{K}^{i+j}  \\
\\
 \left(  \begin{array}{cc}
\pi_{K}^{j} \alpha   & \pi_{K}^{j} \beta   \\
\\
0 & \pi_{K}^{j} \delta  \end{array} \right) \left(  \begin{array}{cc}
\pi_{K}^{s} a_{1} & 0 \\
\\
\pi_{K}^{s} c_{1} & \pi_{K}^{s} a_{1}
\end{array} \right) =  \left(  \begin{array}{cc}
\alpha a_{1} + \beta c_{1}  & \beta a_{1} \\
\\
\delta c_{1} &  a_{1} \delta 
\end{array} \right) \pi_{K}^{j+ s} 
\end{array} \]
so $s=i$ and since $\beta \not= 0$ we have $a = a_{1}$ which implies $c = 0$, which is generally not true. Also $U_{K} \cong K^{*} \times {\mathcal O}_{K}$ so that $H^{1}({\rm Gal}(K/F); U_{K}) = \{ 1 \}$ when $K/F$ is tame and $B_{K} \bigcap U_{K} = K^{*}$ so that $H^{1}({\rm Gal}(K/F); B_{K} \cap U_{K}) = \{ 1 \}$. This vanishing of cohomology descends to $K'$ in $F \leq K' \leq K$, an intermediate Galois extension.

In the preceding discussion $G_{K} \leq GL_{2}K$ stabilises a basic choice of $0$-simplex and the conjugate
\[ \begin{array}{l}
 \left(  \begin{array}{cc}
  \pi_{K} & 0 \\
\\
0 & 1   
\end{array} \right)   G_{K}   \left(  \begin{array}{cc}
  \pi_{K}^{-1} & 0 \\
\\
0 & 1   
\end{array} \right) \\
  \\
  =  \{   \left(  \begin{array}{cc}
  \pi_{K} & 0 \\
\\
0 & 1   
\end{array} \right) \pi_{K}^{i}  \left(  \begin{array}{cc}
a  &  b \\
\\
c &  d
\end{array} \right)    \left(  \begin{array}{cc}
\pi_{K}^{-1}    & 0 \\
\\
0   &  1
\end{array} \right)  \}  \\
\\
  = \{  \left(  \begin{array}{cc}
a \pi_{K}  &  b \pi_{K}  \\
\\
c   &  d   
\end{array} \right)    \left(  \begin{array}{cc}
\pi_{K}^{-1}   & 0 \\
\\
0   &  1
\end{array} \right) \}\\
\\
 = \{  \left(  \begin{array}{cc}
a  &  b \pi_{K} \\
\\
c  \pi_{K}^{-1}  &  d   
\end{array} \right) \}. \\
\\
\end{array} \]
which I shall denote by $G'_{K}$.

 The intersection
\[ \begin{array}{l}
G''_{K} = G_{K}  \bigcap G'_{K} \\
\\
=  K^{*} \cdot GL_{2}{\mathcal O}_{K} \bigcap  \{ K^{*} \cdot \left( \begin{array}{cc}
a & b \pi_{K} \\
\\
c \pi_{K}^{-1} & d
\end{array} \right)  \ | \ a,b,c,d \in {\mathcal O}_{K}, ad-bc \in {\mathcal O}_{K}^{*} 
  \}  \\
  \\
  = \{ \left( \begin{array}{cc}
 \pi_{K}^{i}  a & b \pi_{K}^{i+1}  \\
\\
c \pi_{K}^{i-1} & d  \pi_{K}^{i}
\end{array} \right)  \ | \ a,b,c,d \in {\mathcal O}_{K}, ad-bc \in {\mathcal O}_{K}^{*} , \ c \in \pi_{K} {\mathcal O}_{K}
  \} 
  \end{array}   \]
  so we put 
  \[  \begin{array}{l}
  U''_{K} =  \{ \left( \begin{array}{cc}
 \pi_{K}^{i}  a & 0 \\
\\
c \pi_{K}^{i-1} & a  \pi_{K}^{i}
\end{array} \right)  \ | \ a,c,d \in {\mathcal O}_{K}, ad \in {\mathcal O}_{K}^{*} , \ c \in \pi_{K} {\mathcal O}_{K}
  \} \\
  \\
  B''_{K} =  \{ \left( \begin{array}{cc}
 \pi_{K}^{i}  a & b \pi_{K}^{i+1}  \\
\\
0  & d  \pi_{K}^{i}
\end{array} \right)  \ | \ a,b,d \in {\mathcal O}_{K}, ad \in {\mathcal O}_{K}^{*} 
  \} 
    \end{array}   \]
 The previous discussion shows that $U''_{K}$ selfnormalises in $G''_{K} = U''_{K}  \cdot B''_{K}$ and that the appropriate pair of first Galois cohomology groups both vanish.
 \vspace{5pt}
 \newpage
 Let    
 \[ N_{K}  = \{  \left( \begin{array}{cc}
    1 & \beta \\
    \\
    0 & 1 
    \end{array} \right) \ | \ \beta \in {\mathcal O}_{K}  \}  \]
    and set $B_{K}/N_{K} = T_{K}$ amd similarly with one or two prime superscripts.
 We have ${\rm Gal}(K/F)$ with $F \leq  K$ tame or unramified.
Now we want to produce a (unique) conjugacy class in each of
${\rm Gal}(K/KF) \propto  T_{K}$, ${\rm Gal}(K/KF) \propto  T'_{K}$ and ${\rm Gal}(K/KF) \propto  T''_{K}$.

The three compact open modulo the centre subgroups $G_{K}, G'_{K}$ and $G''_{K}$ are the stabilisers of the three simplices comprising the fundamental domain line-segment (see \cite{Sn18} Chapter Two).

\section{$\pi$ on $GL_{3}$}

In this example the first vertex has stabiliser $G_{K} = K^{*} \cdot GL_{3}{\mathcal O}_{K}$ whose elements look like
\[ G_{K} = \{ \pi_{K}^{i}  X = \pi_{K}^{i}   \left(  \begin{array}{ccc}
a & b &c \\
\\
d & e & f \\
\\
g & h & k 
\end{array} \right)  \ | \ {\rm det}(X)  \in {\mathcal O}_{K}^{*} \}  .  \]

We have a Chevalley subgroup $U_{K} \leq G_{K}$ as in (\cite{LLProb1970} Lemma 1, p.30) of the form
\[    U_{K} = \{ \pi_{K}^{i}  \left(  \begin{array}{ccc}
a & 0&0 \\
\\
d & a & 0 \\
\\
g & h & a 
\end{array} \right) \in G_{K} , \  i \in {\mathbb Z}    \}      \]
and the Borel subgroup
\[    B_{K} = \{   \pi_{K}^{i}  \left(  \begin{array}{ccc}
a & b &c \\
\\
0 & e & f \\
\\
0& 0 & k 
\end{array} \right)  \in G_{K}  , \  i \in {\mathbb Z}    \}      \]
so that $G_{K} = U_{K} B_{K}$ because
\[ \begin{array}{l}
 \pi_{K}^{i}  \left(  \begin{array}{ccc}
a & 0&0 \\
\\
d & a & 0 \\
\\
g & h & a 
\end{array} \right) \pi_{K}^{j}  \left(  \begin{array}{ccc}
a' & b' &c' \\
\\
0 & e' & f' \\
\\
0& 0 & k' 
\end{array} \right)  \\
\\
=  \pi_{K}^{i+j}  \left(  \begin{array}{ccc}
a a' & a b' & a c' \\
\\
d a' &d b' + a e' & d c' +  a f' \\
\\
g a'& g b' + h e' & g c' +  h f' + a k' 
\end{array} \right) .
\end{array} \]

In the preceding discussion $G_{K} \leq GL_{2}K$ stabilises a basic choice of $0$-simplex and the conjugate
\[ \begin{array}{l}
 \left(  \begin{array}{ccc}
  \pi_{K} & 0 & 0  \\
\\
0 & 1 &  0 \\
\\
0 & 0 & 1
\end{array} \right)   G_{K}   \left(  \begin{array}{ccc}
  \pi_{K}^{-1} & 0  & 0 \\
\\
0 & 1 &  0 \\
\\
0 & 0 & 1 
\end{array} \right) \\
  \\
  =  \{   \left(  \begin{array}{ccc}
  \pi_{K} & 0 & 0  \\
\\
0 & 1 &  0 \\
\\
0 & 0 & 1
\end{array} \right)   \pi_{K}^{i}    \left(  \begin{array}{ccc}
a & b &c \\
\\
d & e & f \\
\\
g & h & k 
\end{array} \right)  \left(  \begin{array}{ccc}
  \pi_{K}^{-1} & 0  & 0 \\
\\
0 & 1 &  0 \\
\\
0 & 0 & 1 
\end{array} \right)  \}  \\
\\
  =  \pi_{K}^{i}  \left(  \begin{array}{ccc}
a & \pi_{K} b & \pi_{K}c \\
\\
\pi_{K}^{-1} d & e & f \\
\\
\pi_{K}^{-1} g & h & k 
\end{array} \right)  
\end{array} \]
which I shall denote by $G'_{K}$ and we set $G'''_{K} = G_{K} \bigcap G'_{K}$, which stabiises one of the $1$-simplices of the fundamental $2$-simplex.
\[ \begin{array}{l}
 \left(  \begin{array}{ccc}
1 & 0 & 0  \\
\\
0 &   \pi_{K} &  0 \\
\\
0 & 0 & 1
\end{array} \right)   G_{K}   \left(  \begin{array}{ccc}
  1 & 0  & 0 \\
\\
0 &   \pi_{K}^{-1} &  0 \\
\\
0 & 0 & 1 
\end{array} \right) \\
  \\
  =  \{   \left(  \begin{array}{ccc}
  1 & 0 & 0  \\
\\
0 & \pi_{K} &  0 \\
\\
0 & 0 & 1
\end{array} \right)   \pi_{K}^{i}    \left(  \begin{array}{ccc}
a & b &c \\
\\
d & e & f \\
\\
g & h & k 
\end{array} \right)  \left(  \begin{array}{ccc}
 1 & 0  & 0 \\
\\
0 &  \pi_{K}^{-1} &  0 \\
\\
0 & 0 & 1 
\end{array} \right)  \}  \\
\\
  =  \pi_{K}^{i}  \left(  \begin{array}{ccc}
a & \pi_{K}^{-1}  b & c \\
\\
\pi_{K} d & e & \pi_{K} f \\
\\
g & h \pi_{K}^{-1} & k 
\end{array} \right)  
\end{array} \]
which gives $G''_{K}$ and $G^{iv}_{K} = G_{K} \bigcap G''_{K}$.

\[ \begin{array}{l}
 \left(  \begin{array}{ccc}
  1& 0 & 0  \\
\\
0 & 1 &  0 \\
\\
0 & 0 & \pi_{K} 
\end{array} \right)   G_{K}   \left(  \begin{array}{ccc}
 1 & 0  & 0 \\
\\
0 & 1 &  0 \\
\\
0 & 0 & \pi_{K}^{-1} 
\end{array} \right) \\
  \\
  =  \{   \left(  \begin{array}{ccc}
  1 & 0 & 0  \\
\\
0 & 1 &  0 \\
\\
0 & 0 &   \pi_{K}
\end{array} \right)   \pi_{K}^{i}    \left(  \begin{array}{ccc}
a & b &c \\
\\
d & e & f \\
\\
g & h & k 
\end{array} \right)  \left(  \begin{array}{ccc}
1 & 0  & 0 \\
\\
0 & 1 &  0 \\
\\
0 & 0 &   \pi_{K}^{-1}
\end{array} \right)  \}  \\
\\
  =  \pi_{K}^{i}  \left(  \begin{array}{ccc}
a &  b & \pi_{K}^{-1} c \\
\\
 d & e & \pi_{K}^{-1} f \\
\\
\pi_{K} g & \pi_{K} h & k 
\end{array} \right)  
\end{array} \]
which gives $G^{''}_{K}$ and $G^{v}_{K} = G'_{K} \bigcap G''_{K}$. Finally $G^{vi}_{K} = G_{K} \bigcap G' \bigcap G''_{K}$.
These seven compact open modulo the centre subgroups are the stabilisers of the seven simplices comprising the triangular fundamental domain(see \cite{PG97}  Chapter Seven, flag complexes for $GL_{n}$; particularly p.196 for $GL_{3}$). In general for $GL_{n}$ of a local field the building has $2^{n}-1$ compact modulo the centre stabiliser groups.

\section{Vanishing $H^{1}$-conditions for $GL_{2}$}

{\bf Lang's Theorem and $H_{K}$ for $GL_{2}$}

I am interested in Lang's Theorem for $H_{K}$ when $K/F$ is unramified. This is related to but not equivalent to $H^{1}({\rm Gal}(K/F); H_{K} ) = \{ 1 \}$. I am hoping that if the Lang surjection is onto for $G_{K}, G'_{K}$ and $G''_{K}$ in the $GL_{2}$ example  (and $H_{1} \bigcap H_{2}$) that this will give me in a straightforward cohomological way the results \cite{LLProb1970} on conjugacy classes from the results of \S3.

Recall Lang's Theorem. For $G$ over ${\mathbb F}_{q}$ smooth and ${\rm Frob}_{F_{q}}$ equalling the Frobenius, temporarily denoted by  $\sigma $ for the rest of this paragraph, then the map $x \mapsto x^{-1} \sigma(x)$ is onto. Here is the proof: Consider $f_{a} : x \mapsto x^{-1} a \sigma(x)$ whose differential is $(df_{a})(1) = d( h \cdot (x^{-1}, a, \sigma(x))$  where $h(x,y,z) = xyz$. Hence
\[   (df_{a})(1) = dh_{(1,a,1)} \cdot (-1, 0, d\sigma_{1}) = -1 + d\sigma_{1} = -1 .\]
Now $f_{a}(bx) = f_{f_{a}(b)}(x)$ so $df_{b}$ is a bijection for any $b$. The dimensions of the tangent space and of $G$ at any point are equal so there is a dense open set $U$ containing $b$ and another $V$ containg the image of $f_{1}$. Hence the Lang has a dense open image and is therefore onto.

We know that $H^{i}({\rm Gal}(K/F); {\mathbb F}_{q} ) = \{ 1 \}$ for all $i$ and , since $K/F$ is unramified,
$H^{i}({\rm Gal}(K/F); {\mathcal O}_{ F} ) = \{ 1 \}$  for all $i$. Therefore, by induction 
$H^{i}({\rm Gal}(K/F); {\mathcal P}_{ F}^{n} ) = \{ 1 \}$ for all $i$ and for $n = 0,1, 2 , \ldots $.

For $n \geq 1$ 
\[ 0  \longrightarrow  1 + M_{s} {\mathcal P}_{ F}^{n+1 }   \longrightarrow  1 + M_{s} {\mathcal P}_{ F}^{n } 
\longrightarrow  M_{s} {\mathbb F}_{ q}   \longrightarrow 0 \]
shows that for $n \geq 1$ 
\[ H^{i}({\rm Gal}(K/F);  1 + M_{s} {\mathcal P}_{ F}^{n+1 } )   \stackrel{\cong}{\longrightarrow}
 H^{i}({\rm Gal}(K/F);  1 + M_{s} {\mathcal P}_{ F}^{n } )  . \]
 
 By the Mittage-Leffler condition, for all $i$,
 \[ \begin{array}{l}
 H^{i}({\rm Gal}(K/F); GL_{s}{\mathcal O}_{K})  \\
 \\
\cong    {\rm lim}_{   \stackrel{n}{ \leftarrow} }   H^{i}({\rm Gal}(K/F); GL_{s}{\mathcal O}_{K}/1 + M_{s} {\mathcal P}_{K}^{n} ) \\
\\
= \{ 1 \} 
 \end{array} \]
 when $K/F$ is unramified.
 
 Recall how vanishing $H^{1}({\rm Gal}(K/F) ; -)$ and finite orders can give Lang's Theorem.
 Take the Galois module ${\overline{\mathbb F}}_{q}$ with Frobenius 
 $\sigma \in {\rm Gal}({\overline{\mathbb F}}_{q}/ {\mathbb F}_{q})$. Take $0 \not= x \in {\overline{\mathbb F}}_{q}$
 so that  $x \in {\mathbb F}_{q^{d}}$ and $x^{q^{d}-1} = 1$. Set $e = q^{d}-1$ and consider 
 $x \in ({\mathbb F}_{q^{de}})^{*}$. Then the norm $N(x) \in {\mathbb F}_{q}$ is the product of $e$ copies of 
 $(x \sigma(x) \sigma^{2}(x) \ldots \ \sigma^{d-1}(x))$ which is trivial because $x^{e} =1$ so $N(x)=1$ and since $H^{1}({\rm Gal}( {\mathbb F}_{q^{de}}/ {\mathbb F}_{q}));  ({\mathbb F}_{q^{de}})^{*}) = \{ 1 \}$ there exists 
 $y \in ({\mathbb F}_{q^{de}})^{*}$ such that $y^{-1} \sigma(y) = x$, which is Lang's Theorem.
 
 Since the elements of $GL_{s}{\mathcal O}_{K}/1 + M_{s} {\mathcal P}_{K}^{n}$ have finite order we can prove Lang's Theorem in the same manner for $GL_{s}{\mathcal O}_{K}/1 + M_{s} {\mathcal P}_{K}^{n}$. Then, by Mittag-Leffler, we see that Lang's Theorem holds for  $GL_{s}{\mathcal O}_{F_{unram}}$ where $F_{unram}/F$ is the maximal unramified extension of $F$.
 
 Once one has Lang surjectivity for  $GL_{s}{\mathcal O}_{F_{unram}}$ the observation of Digne-Michel (\cite{Sn18} p.266-68; see also \S3) applies to give a bijection between conjugacy classes in  $GL_{s}{\mathcal O}_{F_{unram}}$ and in ${\rm Gal}( F_{unram}/K) \propto  GL_{s}{\mathcal O}_{F_{unram}})$.
 
 Here is how we get from Lang's Theorem for $GL_{s}{\mathcal O}_{F_{unram}}$ to $H_{K}$ which is the case
 $F_{unram}^{*} \cdot GL_{s}{\mathcal O}_{F_{unram}}$. Suppose we have $\beta^{-1}X$ with $\beta \in F_{unram}^{*} $ and $X \in GL_{s}{\mathcal O}_{F_{unram}}$. Therefore we have $z $ such that $z^{-1} \sigma(z) = \beta$ and $Y$ such that $Y^{-1} \sigma(Y) = X$ so that $\sigma(z^{-1}) z = \beta^{-1}$. Since the scalar matrices are central we have 
 \[  \begin{array}{ll}
  \beta^{-1} X = \sigma(z^{-1}) z Y^{-1} \sigma(Y) 
  & = Y^{-1} \sigma(Y) \sigma(z^{-1}) z \\
  \\
  & =  z^{-1} z  Y^{-1} \sigma(Y) \sigma(z^{-1}) z \\
  \\
  & =  z^{-1} ((Yz)^{-1} \sigma(Y z^{-1})) z  \\
  \\
  & =  (Yz)^{-1} \sigma(Y z^{-1})
 \end{array} \]
 because $z$ is a scalar matrix.

\section{Langlands 2-variable L-functions (\cite{LLProb1970} pp. 29-34)}

{\bf Part A:  An ingenious construction} (\cite{LLProb1970} pp. 29-34)

To define the local $L$-functions, to prove that almost all primes are unramified and to prove that the product over the unramified primes converges for the real part sufficiently large we need some facts from the reduction theory for groups over local fields (see \cite{BT67}). Much progress has been made in that theory but (at the time of  \cite{LLProb1970}) it was still incomplete. Unfortunately the particular facts needs did not seem then to be in the literature. Very little was lost in just assuming them. For the groups under consideration what is needed is easily verified.

Suppose that $K/F$ is an unramified extension of local fields and $G$ is a quasi-split group over $F$ which splits over $K$. Let $B$ be a Borel subgroup of $G$ and $T$ a Cartan subgroup of $B$ both of which are defined over $F$. Let $v$ be the valuation on $K$. It is a homomorphism from $K^{*}$ whose kernel is the group of units. ${\mathcal O}_{K}^{*}$. If $t \in T_{F}$ let $v(t) \in \hat{L}$ be defined by $< \lambda , v(t) > = v( \lambda(t)) $ for all $\lambda \in L$.

Recall from \S2 that $\hat{L} = {\rm Hom}(L, {\mathbb Z})$ 
and from \S3 that $L$ is homomorphisms on $T = T^{n}$ of the form $\lambda \in L$ given by
\[ \lambda : \left(  \begin{array}{ccccc}
t_{1} & 0 & \ldots & \ldots \ & 0 \\
\\
0 & t_{2} &  \ldots & \ldots \ & 0 \\
\\
\ldots & \ldots & \ldots & \ldots & \ldots  \\
\\
0 & \ldots & \ldots & 0 & t_{n} \\
\end{array} \right)   \mapsto  t_{1}^{m_{1}}  t_{2}^{m_{2}}  \ldots  t_{n}^{m_{n}} \] 
with $m_{i} \in {\mathbb Z}$ we write $\lambda = (m_{1}, \ldots , m_{n})$. Thus $L$ is identified with ${\mathbb Z}^{n}$.

If $\sigma \in {\rm Gal}(K/F)$ then
\[ < \lambda , \sigma v(t) > = < \sigma^{-1} \lambda , v(t) > = v( \sigma^{-1}( \lambda (\sigma t))) = v( \lambda(t))  \]
because $\sigma(t) =t$ and $v( \sigma^{-1} a) = v(a)$ for all $a \in K^{*}$. Thus $v$ is a homomorphism from $T_{F}$ into $\hat{M}$, the group of Galois fixed points  $\hat{L}^{ {\rm Gal}(K/F)}$, which is actually a surjection. 
\begin{lemma}{$_{}$}
\label{5.1a}
\begin{em}

There is a Chevalley lattice in the Lie algebra of $G$ whose stabiliser $U_{K}$ is invariant under ${\rm Gal}(K/F)$. $U_{K}$ is self-normalising. Moreover, $G_{K} = B_{K} U_{K} $, 
$H^{1}( {\rm Gal}(K/F) ; U_{K}) = \{1 \}$ and $H^{1}( {\rm Gal}(K/F) ; B_{K} \bigcap U_{K}) = \{1 \}$. If we choose two such Chevalley lattices with stabilisers $U_{K}$ and $U'_{K}$ then $U_{K}$ and $U'_{K}$ are conjugate in $G_{K}$.

By this we mean that if $gU_{K}g^{-1} = U'_{K}$ for $g \in G_{F}$ then $U'_{K}$ has the same four properties
and and $U'_{F} = G_{F} \bigcap U'_{K} = B_{F}U'_{F}$.

\end{em}
\end{lemma} 

{\bf Proof of conjugacy of $U_{F}$ and $U'_{F}$ in $G_{F}$} \vspace{5pt}  

If $g \in G_{K}$ and $\sigma \in {\rm Gal}(K/F)$, let $g^{\sigma} = \sigma^{-1}(g)$. If $g \in G_{F}$ we may write it as $g = bu$ with $b \in B_{K}$ and $u \in U_{K}$. Then $g^{\sigma} = b^{\sigma} u^{\sigma}$  
and $u^{\sigma} u^{-1} = b^{- \sigma} b $. By Lemma \ref{5.1a} there is $v \in B_{K} \bigcap U_{K}$ such that
$u^{\sigma} u^{-1} = b^{- \sigma} b  = v^{\sigma} v^{-1}$. Then $ b' = bv \in B_{F}$, $u' = v^{-1} u \in U_{F} = G_{F} \bigcap U_{K}$ and $g = b' u'$ so that $G_{F} = B_{F} U_{F}$.

If $g U_{K} g^{-1} = U'_{K}$ for $g \in G_{K}$ then $ g^{\sigma} U_{K} g^{- \sigma} = U'_{K}$ so that $g^{- \sigma} g \in U_{K}$ which is its own normaliser. By Lemma \ref{5.1a} there is $u \in U_{K}$ so that $g^{- \sigma} g = u^{\sigma} u^{-1}$. Then $g_{1} = gu \in G_{F}$ and $g_{1}U_{K} g_{1}^{-1} = U'_{K}$ so that 
$U_{F}$ and $U'_{F}$  are conjugate in $G_{F}$.

Let $C_{c}(G_{F}, U_{F})$ be the set of all compactly supported functions from $G_{F}$ to $U_{F}$ such that $f(gu) = f(ug) = f(g)$ for all $u \in U_{F}, g \in G_{F}$\footnote{$U_{F}$ must be compact if $f \not= 0$ and $f$ is compactly supported.}. $C_{c}(G_{F}, U_{F})$ is a convolution algebra called the Hecke algebra. If $N$ is the unipotent radical of $B$ let $dn$ be a Haar measure on $N_{F}$ and let $\frac{d(bnb^{-1})}{dn}  = \delta(b)$ if $b \in B_{F}$. If $\lambda \in \hat{M}$ choose $t \in T_{F}$ such that $v(t) = \lambda$. If $f \in C_{c}(G_{F}, U_{F})$ set
\[ \hat{f}(\lambda) = \delta^{1/2}(t) ( \int_{N_{F} \bigcap U_{F}}  dn)^{-1} \int_{N_{F}} f(tn) dn .\]
The group ${\rm Gal}(K/F)$ acts on $\Omega$. Let $\Omega^{0} = \Omega^{{\rm Gal}(K/F)}$.
 $\Omega^{0} $ acts on $\hat{M}$. Let $\Lambda(M)$  be the complex group algebra of $\hat{M}$ and let 
 $\Lambda^{0}(\hat{M}) = \Lambda(M)^{\Omega^{0}}$.
 \begin{lemma}{\cite{IS63}}
 \label{5.2a}
 \begin{em}
 
 The map $f \mapsto \hat{f}$ is an isomorphism from $C_{c}(G_{F}, U_{F})$ to 
 $\Lambda(M)^{ \Omega^{0} }$.
 \end{em}
 \end{lemma} 
 
 Pages 31 and 32 of \cite{LLProb1970} assures us that different sets of choices of all this structure leads to canonically compatible data - I shall skip the details.
 
 If $\pi$ is an irreducible unitary representation of $G_{F}$ on $H$ whose restriction to $U_{F}$ contains the identity representation then 
 \[ H_{0} = \{ x \in H \ | \ \pi(u)x = x \ {\rm for \ all } \ u \in U_{F} \}  \]
 is one-dimensional. If $f \in C_{c}(G_{F}, U_{F})$ then
 \[   \pi(f) = \int_{G_{F}}  f(g) \pi(g) dg  \]
 maps $H_{0}$ to itself. The representation of $C_{c}(G_{F}, U_{F})$ on $H_{0}$ determines a ring homomorphism $\chi : C_{c}(G_{F}, U_{F}) \longrightarrow {\mathbb C}$ and also 
 $\Lambda^{0}(\hat{M}  )   \longrightarrow  {\mathbb C}$. To define the local L-functions we study such homomorphisms. First of all observe that if $\chi$ is associated to a unitary representation then 
 \[   | \chi(f) | \leq  \int_{G_{F}}  |f(g)| dg . \]
 since  $\Lambda(\hat{M})  $ is finitely generated over  $\Lambda^{0}(\hat{M})  $ any homomorphism like $\chi$ can be extended to a ring homomorphism of the form  $\Lambda(\hat{M}  )   \longrightarrow  {\mathbb C}$.
which will necessarily be of the form $ \Sigma \hat{f}(\lambda) \lambda \mapsto   \Sigma \hat{f}(\lambda) \lambda(t)$ for some $t \in \hat{T}$. Conversely given $T$ the formula above determined a ring homomorphism $ \chi_{t}: \Lambda^{0}(\hat{M}  )   \longrightarrow  {\mathbb C}$. We shall show that $\chi_{t_{1}} = \chi_{t_{2}}$ if and only if $t_{1} \times \sigma_{F} $ and $ t_{2} \times \sigma_{F}$, where $\sigma_{F}$ is the Frobenius substitution, are conjugate in $\hat{G}_{F}$. We abbreviate $t \times \sigma$ to $t \sigma$. Every semi-simple conjugacy class of  $\hat{G}_{F}$ is conjugate to some $ t \sigma_{F}$ with $t \in \hat{T}$ \cite{FG39}.

Thus there is a bijection between complex-valued ring homomorphisms of the Hecke algebra and semi-simple conjugacy classes in $\hat{G}_{F}$ whose projection to ${\rm Gal}(K/F)$ is $\sigma_{F}$.

If $\rho$ is a complex analytic representation of $\hat{G}_{F}$ and 
$ \chi_{t}: \Lambda^{0}(\hat{M}  )   \longrightarrow  {\mathbb C}$ is the homomorphism associated to $\pi$ we define the local L-function to be
\[   L(s, \rho, \pi) = \frac{1}{{\rm det}( 1 - \rho( t \sigma_{F}) | \pi_{F}|^{s})}  .  \]

 \begin{example}
 \label{5.2a}
 \begin{em}

The map $f \mapsto \hat{f}$ is an isomorphism from  $C_{c}(G_{F}, U_{F})$ to $\Lambda^{0}(M)$.
In this example $U_{K} = K^{*} = B_{K}$.
\end{em}
\end{example}

{\bf   Part B: Left actions}

Let $K$ be an unramified extension of $F$.

In \cite{Sn20} my convention for the structure of the Galois semi-direct product in Lemma \ref{5.1a} is
$(\sigma, g) (\sigma'.g') = (\sigma \sigma' ,  g \sigma(g'))$ . The monomial resolution for $\tilde{\pi}$ on ${\rm Gal}(K/F) \propto GL_{n}(K)$, the base change of $\pi$, is a sum indexed by stabilisers $G_{K}$ of the individual simplices in a fundamental domain of $GL_{n}F$ acting on the Bruhat-Tits building \footnote{If $G_{K}$  is a simplex stabiliser for  $\pi$ then ${\rm Gal}(K/F) \propto G_{K}$ is a set of simplex stabilisers for the base-change, Galois-permuted in the obvious manner.}. The function I attach to $\rho$, $\pi$ is the Euler product of functions, one for each fundamental simplex stabiliser, which Langlands gives canonically for each piece of $\pi$ associated to a "special Langlands line".

From the data, monomial resolution of $\tilde{\pi}$, simplex stabiliser, special Langlands line, Langlands gives us, in $Gal(K/F) \propto GL_n(K)$, $(Frob_{P},t)$ canonical, with Galois coordinate presumably equal to $Frob_{P}$. $t$ is an $F$-point of a torus $T$ dependent on the data.The Galois norm of $(Frob_{P},t)$ will be equal to $(Frob_{K},t^{[K:F]})$ which contributes
\[\frac{1}{{\rm det} ( 1 - \rho( Frob_{K},t^{[K:F]}) | \pi_{K}|^{s})}, \]

I believe. 

\vspace{5pt}
{\bf Question}: What are the analogous Galois norms from the Euler characteristic of the remaining data from the monomial resolution?

 \end{document}